\documentclass{amsart}
\usepackage[margin=1.5in]{geometry}
\usepackage{graphicx}
\usepackage[svgnames]{xcolor}
\usepackage[colorlinks,urlcolor=NavyBlue,citecolor=Maroon,linkcolor=DarkSlateGrey]{hyperref}
\usepackage{amsthm}
\theoremstyle{plain}

\theoremstyle{definition}

\usepackage{amssymb}
\usepackage{amsmath}
\usepackage{tikz}
\usetikzlibrary{matrix,arrows,positioning,automata}

\newcommand{\T}[1]{\textbf{\texttt{#1}}}

\newcommand{\cI}{\mathcal{I}}

 \title[The Attractor-Cycle Notation]{The Attractor-Cycle Notation for Finite Transformations}
 \author{Attila Egri-Nagy$^1$}
 \address{$^1$Akita International University, Japan}
 \email{egri-nagy@aiu.ac.jp}
 \author{Chrystopher L. Nehaniv$^2$}
 \address{$^2$University of Waterloo, Canada}
 \email{cnehaniv@uwaterloo.ca}

\begin{document}

\maketitle

\begin{abstract}
  We describe a new notation for finite transformations.
This attractor-cycle notation extends the orbit-cycle notation for permutations and  builds upon existing transformation notations.
How the basins of attraction of a finite transformation flow into permuted orbit cycles is visible from the notation.
It gives insight into the structure of transformations and reduces the length of expressions without increasing the number of types of symbols.
\end{abstract}
\section{Introduction}

\begin{figure}
\begin{center}
  \tikzset{->,>=triangle 45,auto,node distance=1.5cm}
\begin{tikzpicture}[scale=0.9, every node/.style={scale=0.9}]
\tikzstyle{every state}=[minimum size=1pt]
  \node[state]  (q_1)  {$1$};
  \node[state]  (q_2) at(.75,-1) {$2$};
  \node[state]  (q_3) [right of=q_1] {$3$};
  \node[state]  (q_4) [below of=q_2] {$4$};
  \node[state]  (q_5) [below of=q_3] {$5$};
  \node[state]  (q_6) [below of=q_4] {$6$};

  \node[state]  (q_7) at(2.5,0) {$7$};
  \node[state]  (q_8) [right of=q_7] {$8$};
  \node[state]  (q_9) at(3,-1.5) {$9$};
  \node[state]  (q_10) [below of=q_9] {$10$};
  \node[state]  (q_11) at(4.5,-1.5) {$11$};
  \node[state]  (q_12) at(6,-3) {$12$};
  \node[state]  (q_13) [above of=q_12] {$13$};
  \node[state]  (q_14) [above of=q_13] {$14$};
  \node[state]  (q_15) [right of=q_12] {$15$};
  \node[state]  (q_16) [above=1cm of q_15] {$16$};
  \node[state]  (q_17) [right of=q_16] {$17$};

  \path (q_16) edge  [bend left] (q_17);
  \path (q_17) edge  [bend left] (q_16);
  \path (q_1) edge  (q_2);
  \path (q_3) edge  (q_2);
  \path (q_2) edge  (q_4);
  \path (q_5) edge  (q_4);
  \path (q_4) edge  (q_6);

  \path (q_10) edge  [bend left](q_11); 
  \path (q_11) edge  [bend left](q_12);
  \path (q_12) edge  [bend left] (q_10);
  \path (q_14) edge  (q_13);
  \path (q_13) edge   (q_12);
  \path (q_7) edge   (q_9);
  \path (q_8) edge   (q_9);
  \path (q_9) edge   (q_10);
  \path (q_6) edge [loop,in=225,out=315,looseness=9] (q_6);
  \path (q_15) edge [loop]  (q_15);
\end{tikzpicture}

\caption{A simple discrete dynamical system. A transformation $g:X \rightarrow
  X$ with four components, i.e.\ basins of attraction (generalized cycles), is
  visualized as a digraph with arrows $(x,g(x))$ for each $x\in
  X=\{1,\ldots,17\}$. In our canonical form of the attractor-cycle notation, $g$
  is denoted \T{[[[1|3,2]|5,4],6]([[7|8,9],10],11,[14,13,12])(16,17)}. Note:
  Singleton components are not written, just as in orbit-cycle notation for
  permutations. Full, formal details for  the notation are given in Section \ref{acnotation}.}
\label{example1}
\end{center}
\end{figure}
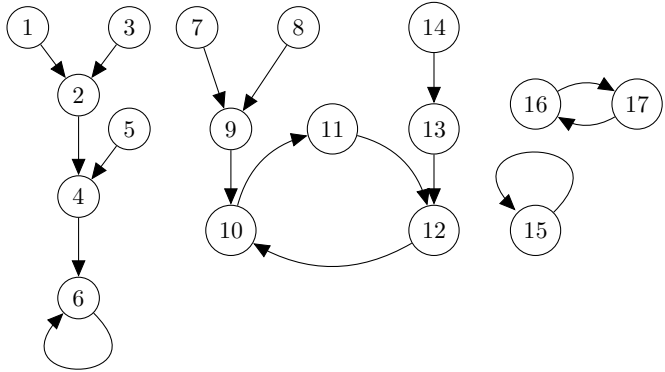

What is the correct generalization of the orbit-cycle notation for permutations to the case of arbitrary total functions $f:X\rightarrow X$ on a finite set $X$?
The discrete dynamical system that $f$ gives on $X$ can be visualized as digraph with nodes $x\in X$ and directed edges $(x,f(x))$.  
Iterating $f$ maps a point to $x$ to $f(x)$, to $f(f(x))$, and so on, until eventually some $f^m(x)=f^{m+k}(x)$.
Taking the least   $m\geq 0$ for which this happens, and the least positive $k$ for that $m$, shows that $x$ eventually must enter a (possibly degenerate) periodic orbit from which it never leaves.
Points $x$ and $x'$ eventually entering the same periodic cycle are said to be in  the same {\it basin of attraction} or {\it generalized cycle}.
Drawn as digraphs, transformations may have several such disjoint basins of
attraction, each consisting of  a cycle of points with incoming trees (connected
acyclic subgraphs). Figure \ref{example1} gives an example with four generalized
cycle components in our notation, which is formally introduced in
Section~\ref{acnotation} with further examples.
For a permutation, this digraph is a set of disjoint cycles.
However, for a transformation, the points in a cycle can have incoming edges from outside the cycle, so a notation for transformations has to deal with these tree structures.
The trees, which could also be degenerate,  are directed toward the cycle, and there may be a tree of points coming into any given point of the cycle.

We call each basin of attraction associated to $f$, a generalized cycle or  a \emph{component}, since it a connected component of the digraph of $f$.  We may restrict our attention to a single component only, since we can write $f$ by concatenating our notation for what $f$ does on each of its components.

The basins of attractions have been studied in the context of dynamical systems as they describe the global dynamics of a temporally evolving system
\cite{wuensche1992global}.
It is important to note that this is the only dynamics in systems where the only action is the unidirectional passage of discrete time.\footnote{This is in contrast to more general discrete dynamical systems which are generated by several such transformations on a finite set $X$, as studied in transformation semigroup theory.}

Various previous notations have been developed for representing such transformations on a set (Section~\ref{othernotations}).
The aim of these notations is to give useful information about the transformation without drawing the corresponding digraph.
Readability is not an easily measurable quantity, but length and the number of distinct symbols used are influencing factors.
The real importance of an efficient notation for transformations lies in the growing  use of computer algebra systems in finite semigroup theory, as well as in mathematical calculations with transformations. 

\section{Existing Notations}\label{othernotations}

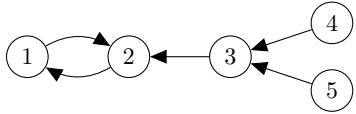
\begin{figure}
\begin{center}    
\tikzset{->,>=triangle 45,auto,node distance=1.5cm}
\begin{tikzpicture}[scale=0.9, every node/.style={scale=0.9}]
\tikzstyle{every state}=[minimum size=1pt]
  \node[state]  (q_1)  {$1$};
  \node[state]  (q_2) [right of=q_1] {$2$};
  \node[state]  (q_3) [right of=q_2] {$3$};
  \node[state]  (q_4) at(4.5,.5) {$4$};
  \node[state]  (q_5) at(4.5,-.5) {$5$};

  \path (q_1) edge  [bend left] (q_2);
  \path (q_2) edge  [bend left] (q_1);

  \path (q_3) edge  (q_2);
  \path (q_4) edge  (q_3);
  \path (q_5) edge  (q_3);

\end{tikzpicture}
 \end{center}
  
\caption{An  irreversible transformation on 5 points with just one basin of attraction.  Here a non-degenerate tree collapses in two iterations into a  cyclic permutation $(1,2)$.}
  \label{runningexample}
\end{figure}

Several ways for generalizing cyclic notation have been developed based on the particular purposes of each situation in which they were defined.\footnote{Somewhat ironically, we also developed yet another version (called compact notation), before settling on the attractor-cycle notation presented here.} Here we give a quick summary of previous suggestions. Throughout this paper, writing $x \mapsto f(x)$ denotes as usual that $f$ maps $x$ to $f(x)$.  Our  running example,  whose digraph is visualized in Figure~\ref{runningexample}, with  just one nontrivial component (basin of attraction) will be the transformation typically written in so-called Cayley notation as
\begin{equation*}   
\left( \begin{array}{ccccc}
1 & 2 & 3 & 4 & 5 \\
2 & 1 & 2 & 3 & 3 
\end{array} \right).
\end{equation*}
This common notation denotes the function on the set $\{1,2,3,4,5\}$ taking an element  $x$ in the first row to the element written $f(x)$ immediately below it on the second row. 

\noindent
Caveat: As it often happens in mathematics, the  notations surveyed in this section associate different meanings to the same symbols.

\subsection{Path Notation for Partial Symmetries and Partial Transformations} In addition to parentheses \T{(},\T{)} for permutation orbits, path notation introduced by S.~Lipscomb~\cite{lipscomb1996symmetric,lipscombbookreview1997} uses square bracket \T{]} following elements with no images in partial permutations.
For partial transformations the paths connect into other paths and ultimately into the cycles.
This is denoted by the symbol \T{>}, may be used as a visual indication of being funneled into a cycle.
The example written in path notation is
\begin{center}\T{(1,2)(4,3,2>(5,3,2>},\end{center}
which redundantly represents the paths going into point 2.

\subsection{Factorization Notation} In order to avoid path redundancy, G.~Ayik, H.~Ayik, and J.M.~Howie~\cite{factorisation_notation_2005_Howie} introduced a different notation.
Instead of decomposing the transformation into paths and cycles, this notation decomposes it into a very particular type of generalized cycle, i.e.\  into a product of  transformations given by the trajectory of a single point.
The example written in path notation is
\begin{center}\T{[4,3,2,1|2][5,3|3]}$=\left(\begin{smallmatrix}1&2&3&4&5\\2&1&2&3&3\end{smallmatrix}\right)$,\end{center}
where the element after the \T{|} shows where the path connects to itself to form a cycle.
This notation is excellent for describing factorizations in the full transformation semigroups, but it introduces maps that might appear to move things in ways that are not present in the original transformation, since branches connecting to an already described path are terminated with an artificial self-loop.
For instance, $3\mapsto 3$.

Also, the decomposition is not unique.
\subsection{Linear Total Transformations}
In \cite{ClassicalTransSemigroups2009} a new notation was introduced by O.~Ganyushkin and V.~Mazorchuk, aiming to provide a natural extension of the cyclic notation of permutations. Instead of decomposing trees into paths, linear notation describes the trees explicitly.
Trees with one level of branching are denoted by \T{[}preimages of root  \T{;} root \T{]},
where preimages are separated by commas and all map to the root.
If a preimage  element also has incoming edges from other points, then  the same square bracket structure is applied again recursively to specify the tree leading to that element.  
Although called `linear'   as  a one-line notation, the notation describes recursive tree structure.

Parentheses are used to indicate the existence of a nontrivial permutation of the roots of the trees (which may include degenerate trees consisting of a single point).
Basically the linear notation is the usual cyclic notation used for permutations but the elements in the cycle describe their incoming tree information. The running example is written as
\begin{center}\T{([[4,5;3];2],1)}.\end{center}

One of the most useful features of linear notation is that when restricted to permutations it is identical to orbit-cycle form. Also, by looking for parentheses we can easily spot the existence of nontrivial permutations even in large examples.
The only drawback is that describing a transformation given by following a simple path, 
a ``line'' of maps, requires many square brackets.
For instance,
\begin{center}
\tikzset{->,>=triangle 45,auto,node distance=1.5cm}
\begin{tikzpicture}[scale=0.9, every node/.style={scale=0.9}]
\tikzstyle{every state}=[minimum size=1pt]
  \node[state]  (q_1)  {$1$};
  \node[state]  (q_2) [right of=q_1] {$2$};
  \node[state]  (q_3) [right of=q_2] {$3$};
  \node[state]  (q_4) [right of=q_3] {$4$};
  \node[state]  (q_5) [right of=q_4] {$5$};

  \path (q_1) edge  (q_2);
  \path (q_2) edge  (q_3);
  \path (q_3) edge  (q_4);
  \path (q_4) edge  (q_5);

\end{tikzpicture}
\end{center}
becomes \T{[[[[1;2];3];4];5]}.
With the attractor-cycle notation we aim to alleviate this particular problem.

\subsection{Specialized Notation for Fast Iterations}
For special purposes one can develop specialized notations.
For example in \cite{FastIterations2007}, for the purpose of the efficient
calculation of $f^m(x)$, the example would be first decomposed into components
\T{((1,2),(4,3),(5))}, then the orbits concatenated \T{(1,2,4,3,5)}  in order to
assign indices to points. This of course lacks the information on how to put together the orbits, therefore an auxiliary sequence is needed to encode where (which position) the last elements of the orbit components connect to: $f(2)$ is in position 1, $f(3)$ in 2 and $f(5)$ in 4. So the full information is encoded in these strings:
\begin{center}
\T{((1,2),(4,3),(5))}, \T{(1,2,4,3,5)}, \T{(1,2,4)}.
\end{center}
This illustrates that efficient algorithmic representations do not necessarily provide human readable formats.

\section{The Attractor-Cycle Notation}\label{acnotation}

Informally, the attractor-cycle notation makes use of the following conventions and symbols.  Formal semantics of the notation are in the next subsection.   Here we illustrate the concepts for points:
\begin{description}
    \item[\textbf{conveyor belt}] Left-to-right comma-separated points. For example, \T{1,2,3} reads as $1\mapsto 2$, $2\mapsto 3$.

    \item[\textbf{cycle}] Putting parentheses around a conveyor belt connects the last point to the first, i.e., \T{(1,}\ldots\T{,n)} adds the map $n\mapsto 1$.

    \item[\textbf{tree operator}] Square brackets group points and do not specify the image of the last point. Based on context, this is resolved into one of 3 possibilities.
    \begin{enumerate}
        \item \T{[}\ldots\T{n],k}\ldots defines the map $n\mapsto k$, connecting into a conveyor belt;
        \item \T{(m,\ldots[}\ldots\T{n])} defines \T{n} $\mapsto$ \T{m}, closing a cycle;
        \item \T{[}\ldots\T{n]} corresponds to $n\mapsto n$, defining a loop edge. 
    \end{enumerate} In short, \T{[\ldots n]} is a tree ``flowing'' into point \T{n}, and the context of \T{n} defines what happen to \T{n} itself.

    \item[\textbf{parallel paths}] A vertical bar \T{|} suggestive of \textbf{logical OR} is written to separate two or more points (respectively, subtrees) that go to the same point. For instance, \T{1|2|3,4} means $1 \mapsto 4$, $2 \mapsto 4$, and $3 \mapsto 4$.
\end{description}

\subsection{Examples}
Any permutation in the attractor-cycle notation is the same as the cyclic notation.
In particular, the identity is simply \T{()}.

A constant map to $n$ is written as \T{[1|2|\ldots|n-1,n]}, while a simple path trajectory from 1 to n is denoted by \T{[1,2,\ldots,n]}. This compact notation for a path is also better matching the cycle notation for permutations.

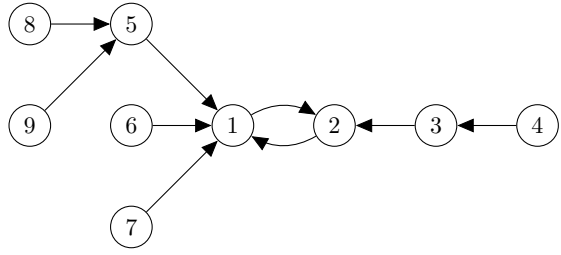
\begin{figure}
    \begin{center}
        \tikzset{->,>=triangle 45,auto,node distance=1.5cm}
\begin{tikzpicture}[scale=0.9, every node/.style={scale=0.9}]
\tikzstyle{every state}=[minimum size=1pt]
  \node[state]  (q_1)  {$1$};
  \node[state]  (q_2) [right of=q_1] {$2$};
  \node[state]  (q_3) [right of=q_2] {$3$};
  \node[state]  (q_4) [right of=q_3] {$4$};
  \node[state]  (q_6) [left of=q_1] {$6$};
  \node[state]  (q_5) [above of=q_6] {$5$};
  \node[state]  (q_7) [below of=q_6] {$7$};
  \node[state]  (q_8) [left of=q_5] {$8$};
  \node[state]  (q_9) [below of=q_8] {$9$};

  \path (q_1) edge  [bend left] (q_2);
  \path (q_2) edge  [bend left] (q_1);
  \path (q_3) edge  (q_2);
  \path (q_4) edge  (q_3);
  \path (q_5) edge  (q_1);
  \path (q_6) edge  (q_1);
  \path (q_7) edge  (q_1);
  \path (q_8) edge  (q_5);
  \path (q_9) edge  (q_5);

\end{tikzpicture}        
    \end{center}
\caption{Digraph for an example transformation
with a `conveyor belt' trajectory $\T{[4,3,2]}$ and a branching tree $\T{[[8|9,5]|6|7,1]}$ collapsing into a transposition $\T{(1,2)}$.  In attractor-cycle notation the whole transformation is denoted $\T{([[8|9,5]|6|7,1],[4,3,2])}$
\label{example3}}
\end{figure}

The  transformation whose digraph is in Figure~\ref{example3} on 9 points is written as
$$\T{([[8|9,5]|6|7,1],[4,3,2])}.$$
The parentheses immediately tell that this transformation has a non-trivial permutation component. Collapsing the tree operators to their roots, we can see that the cycle is \T{(1,2)}.

\vspace{-1em}
\subsection{Syntax} The following context-free grammar defines the language of attractor-cycle notations. The terminal symbols are \T{[},\T{]},\T{(},\T{)},\T{,},\T{|} and the symbols for the $n$ points. The nonterminal symbols are $C$ for components, $F$ for in-flows, $T$ for trees and $P$ for points.
Note the difference between \T{|} and $|$.
The bold symbol \T{|} is used as part of the notational grammar defined here, while
$|$ is part of the metasyntax (Backus-Naur Form (BNF) for context-free grammars) for describing the structure of the new notation. Similarly, parentheses, i.e., $($, $)$, and also $^+$, are part of the metasyntax, while  bold parentheses, i.e., $\T{(}$ and $\T{)},$ are terminal symbols of the attractor-cycle notation.
\begin{center}
    \begin{align}
    S &\rightarrow  C^+ \mid \T{()}\label{rul:comps}\\
    C &\rightarrow \T{(} (T\T{,})^+T\T{)}\mid F \label{rul:comp}\\
    F &\rightarrow \T{[} (T\T{,})^+P\T{]}\mid \T{[} (T\T{|})^+T\T{,}  P \T{]}\label{rul:in-flow}\\
    T &\rightarrow  F\mid P\label{rul:tree}\\
    P &\rightarrow  \T{1}\mid\T{2}\mid\T{3}\mid\ldots\mid\T{n}\label{rul:point}
    \end{align}
  \end{center}
Rule (\ref{rul:comps}) states that the attractor-cycle notation denotes a positive number nontrivial components (generalized cycles, i.e.\ basins of attraction), or we can have just an empty cycle for the identity transformation.
Rule (\ref{rul:comp}) says that a component can be a in-flow (where something flows into a point) or a permutation cycle (with incoming trees to its members).
Rule (\ref{rul:in-flow}) describes the two cases of an in-flow: either a tree with multiple branches leading to the same root  (parallel branches) or a  path (conveyor belt).
Rule (\ref{rul:tree}) defines a tree to be an in-flow. It allows a point to be a tree.
Rule (\ref{rul:point}) specifies the points. 

An important additional constraint in the notation is that each point can only occur at most once.\footnote{ 
As the intersection of a regular language with a context-free one is still context-free \cite{CFintersectREG}, we may constrain the strings generated by the above context-free grammar to have no repeated points (i.e., by intersecting with the regular language of all strings $w$ over all the symbols we use such that $w$ does not contain any repeated point symbol from among $\T{1},\T{2},\T{3}, \ldots, \T{n}$ ) and obtain a context-free language for the attractor-cycle notation without any repeated points.}

\subsection{Semantics}
We define the semantics by recursively interpreting the valid attractor-cycle notation words as a collection of individual maps of points. We write such a map $p\mapsto q$ as a pair $(p,q)$. We will put the total function $f$ together from these pieces.  Using the parse tree of a well-formed word $w$ in the notation, we denote its interpretation by $\cI(w)$. This determines a unique function $f$ from $\{1,\ldots,n\}$ to itself, by identifying $f$ with the set of pairs $\cI(w)$.  If $w=\T{()}$, then
$\cI(w)$ is the identity transformation.  Otherwise $w$ is derived using $S\rightarrow C_1\  \ldots C_k$, and we define
\begin{equation}\label{TopLevel}
\cI(S):=\bigcup_{i=1}^k\cI(C_i)\cup\{(p,p)\mid f(p)\text{ is not defined in any }\cI(C_i)\}.
\end{equation}
We label subtrees to define an auxiliary function $r$ that gives the root of a tree, so we let
$$r(p)=r(\T{[}T_1\T{|}\ldots\T{|}T_k\T{,}p\T{]})=p \in \{1,\ldots,n\}.$$
$$\mbox{ and } r(\T{[}T_1\T{,}\ldots\T{,}T_k\T{]})=r(T_k).$$ 
Then the interpretation of a component derived from a nonterminal symbol $C$ is:
\begin{equation}
\cI(C):=\begin{cases}
\cI(F)&\mbox{if }C\rightarrow F\\
 & \\
\bigcup_{i=1}^k\cI(T_i)&\mbox{if }C\rightarrow\T{(}T_1\T{,}\ldots\T{,}T_k\T{)}\\
\ \cup\left\{(r(T_i),r(T_{i+1})):1\leq i<k\right\}\\
\ \cup\left\{(r(T_k),r(T_1)) \right\}
\end{cases}
\end{equation}

\noindent A tree gives a nonempty set of pairs, if it is 
 a nontrivial in-flow. A point (degenerate tree) does not yield a map. 
\begin{equation}
\cI(T):=\begin{cases}
\cI(F)&\mbox{if } T\rightarrow F\\
\varnothing &\mbox{if }T\rightarrow P
\end{cases}
\end{equation}
A nontrivial tree is an in-flow and it gives a set of pairs as follows:
\begin{equation}
\cI(F):=\begin{cases}
\bigcup_{i=1}^k\cI(T_i)\cup\left\{(r(T_i),p):1\leq i\leq k \right\}&\mbox{ if } F\rightarrow \T{[}T_1\T{|}\ldots\T{|}T_k\T{,}p\T{]}\\
 & \\
\bigcup_{i=1}^k\cI(T_i) &\mbox{ if } F\rightarrow \T{[}T_1\T{,}\ldots\T{,}T_k\T{]}\\
\ \cup\left\{(r(T_i),r(T_{i+1})):1\leq i<k\right\}
\end{cases}
\end{equation}
The notation $w$ thus clearly determines a unique well-defined transformation $f=\cI(w)$ on $\{1,\ldots,n\}$ to itself, since each element $p$ appears at most once in $w$. The transformation is total (i.e, not partial, but fully defined) by Equation~\ref{TopLevel}. Moreover, every $f$ can be written in this attractor-cycle notation in a canonical form.

\subsection{Canonical Form}

The price to pay for the short length is the loss of uniqueness. Both $\T{[1,2,[3|4,6]]}$ and $\T{[[1,2]|3|4,6]}$ denote the same transformation. One can simply choose between a conveyor belt  (comma-separated list of elements in a trajectory) or the parallel branches (using $\T{|}$).
However, a simple recursive algorithm that starts from the point(s) of each component's cycle can produce a canonical form.
All we need to do is to examine the cardinality of the preimage set of each cycle element  from outside the cycle, i.e.\ the number of incoming arrows.
\begin{enumerate}
\item If there is no incoming arrow, then we have a leaf of the tree and
  recursion stops here.
\item If there is only one preimage, then a conveyor belt is built by traversing
  the tree as long as there is only one preimage. Then recursion is done on the
  first element of the conveyor belt built so far.
\item In case there is more than one element in the preimage then parallel branches need to be used with recursion on all parallel elements.
\end{enumerate}

For each transformation $f$ on $n$ points one can now obtain a completely
canonical expression in this notation by additionally requiring  the $C$'s to
appear in order according to their least elements, and that the trees separated by a $\T{|}$ are ordered by their least element, and the cycles start with their least point (or tree root) first.

For instance, here are the conjugacy class representatives of the full transformation semigroup $T_4$ on four points  in canonical form:\\

\T{[1|2|3,4]} 

\T{[[1,2]|3,4]} 

\T{[1|2,3]} 

\T{[[1|2,3],4]} 

\T{[1,2,3,4]} 

\T{[1,2,3]} 

\T{[1,2][3,4]}

\T{[1,2]} 

\T{[1,2](3,4)} 

\T{()} 

\T{(1,2)} 

\T{([1,2],3)} 

\T{(1,2,3)} 

\T{([1|2,3],4)} 
 
\T{([1,2],[3,4])} 

\T{([1,2,3],4)} 

\T{(1,2)(3,4)} 

\T{([1,2],3,4)} 

\T{(1,2,3,4)}\\

Looking at the list, it is easy to spot the existence of nontrivial cycles and idempotents. Idempotents, other than the identity \T{()}, are exactly those transformations given by concatenating a  number of conveyor belts of length 2,  $\T{[}x,y\T{]}$, and single-level parallel branches $\T{[} x_1\T{|}...\T{|}x_k \T{,} x_{k+1}\T{]}$. 
Also, a common feature of all the notations discussed here is that conjugation by a permutation is just  relabelling of points according to the permutation, e.g. $\T{(1,2,3)}^{-1} \T{([1,2],3,4)}\,\T{(1,2,3)}= \T{([2,3],1,4)}$.

\section{Computer Algebra Implementation}
We implemented the attractor-cycle notation in our semigroup decomposition
computer algebra package \texttt{SgpDec} \cite{SgpDec2014,SgpDec} available for
the \texttt{GAP} system \cite{GAP4}.

To verify the correctness, we take a transformation, construct the
corresponding notation string, then parse this string back to a transformation,
and finally we check whether the input is identical to the output.
We conducted this test systematically for all transformations on up to 8 points
($8^8$ transformations), and then for random transformations on up to $2^{20}$ points
(more than one million).

These tests also indicate the practical usability. Parsing context-free grammars
are known to be efficient (polynomial \cite{1967parsing}), and constructing expressions in the notation is
essentially the reverse process of the recursive parsing algorithm.
\section{Conclusion}
By the spreading use of computer algebra systems for investigating transformations and discrete dynamical systems, efficient notation has become a necessity.
The orbit-cycle has been described in numerous group theory textbooks.
Transformations can have the same status as well.

In the attractor-cycle notation we tried to blend the best features of previous notations and also considered the computational experience to derive what we have found to be a useful and readable notation, giving insight into the structure of transformations.

\section*{Acknowledgment}
This work reported in this article was funded in part by the EU project BIOMICS, contract number CNECT-318202, and by the Natural Sciences and Engineering Research Council of Canada (NSERC), funding reference number RGPIN-2019-04669. 
This support is gratefully acknowledged.

\bibliographystyle{plain}
\bibliography{transfnot,../coords}

\end{document}